\documentclass[12pt]{article}
\usepackage{amsfonts}
\usepackage{amssymb}
\usepackage[dvips]{graphics}
\newtheorem{theorem}{Theorem}[section]

\newcommand{\rem}[1]{{\bf Remark:}~}

\def\QED{{\hspace*{\fill}{\vrule height .5ex width 1ex }\quad} 
    \vskip 0pt plus20pt}
\newcommand{\be}{\begin{equation}}
\newcommand{\ee}{\end{equation}}
\newcommand{\bea}{\begin{eqnarray}}
\newcommand{\eea}{\end{eqnarray}}
\newcommand{\beann}{\begin{eqnarray*}}
\newcommand{\eeann}{\end{eqnarray*}}
\newcommand{\bi}{\begin{itemize}}
\newcommand{\ei}{\end{itemize}}
\DeclareMathAlphabet{\mathol}{OT1}{cmr}{l}{ol}

\begin{document}
%{\baselineskip=10pt \thispagestyle{empty} {{\small Preprint UC Davis Math
%1999-24}
%{\tt cond-mat/0111242n} and {\tt mp\_arc 99-nnn} 

\vspace{20pt}

\title{\vspace*{-.35in}
}
\author{{\bf Oscar Bolina}\\[6pt]
%EndAName
Instituto de F\'{\i}sica\\
Universidade de S\~ao Paulo\\
Caixa Postal 66318 \\
S\~ao Paulo 05315-970 Brasil  \\
{\bf E-mail:} bolina@if.usp.br
}
\title{\vspace{-1in}
{
\bf 
The gambler's ruin problem in path representation form  
}}
\date{}
\maketitle
\begin{abstract}
\noindent 
We analyze the one-dimensional random walk of a particle on 
the right-half real line. The particles starts at $x=k$, 
for $k >0$, and ends up at the origin. We solve for the 
probabilities of absorption at the origin by means of a
geometric representation of this random walk in terms 
of paths on a two-dimensional lattice.
\end{abstract}
\section{Introduction}
\noindent 
We consider the classical one-dimensional random walk of a particle
on the right-half real line. We assume that the particle is 
initially at position $x=k$, $k > 0$, and moves to the right with 
probability {\it p} or to the left with probability $1-p$. We 
consider that the particle is absorbed at the origin without 
fixing the number of steps needed to get there. We calculate 
the probability $P(x=k)$ that the particles end up at the origin, 
given that it starts at $x=k$, by means of a geometric 
representation of this random walk in terms of paths on a 
two-dimensional lattice.
\newline
It is well-known that the probability $P(x=k)$ is the $k^{th}$ power 
of $P(x=1)$:
\be\label{xi}
P(x=k)=\left ( P(x=1) \right )^{k}.
\ee
Because of this, we are not really going to {\it prove} the results 
that lead to (\ref{xi}) (the proofs are quite obvious in our
context). Instead, we prefer to use our geometric 
representation as a "visual" proof of these results.
\newpage
\section{The probabilities of absorption}
\begin{figure}
\begin{center}
\resizebox{!}{8 truecm}{\includegraphics{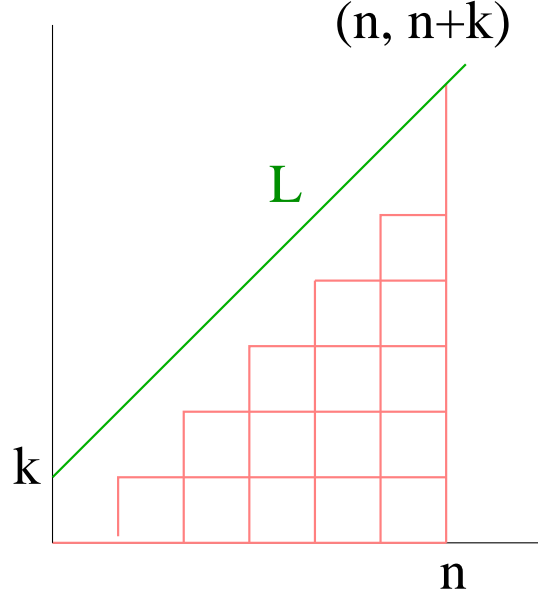}}
\vskip .2 cm
\parbox{13truecm}{\caption{\baselineskip=16 pt\small \label{fig1}
Paths on the lattice, from the origin to the point $(n, n+k)$ on the line
{\it L}, represent the number of steps to the right (n), and to 
the left (n+k) that the particle has taken in its random walk 
from position $x=k$ on the real axis until it is absorbed at the origin. }
}
\end{center}
\end{figure}

A representation of this random walk in terms 
of path on a two-dimensional lattice is given in Fig. \ref{fig1} 
for an arbitrary initial position $x=k$ ($k > 0$) of the 
particle. In this representation, each horizontal lattice 
bond represents a step of the particle to the right, 
while each vertical bond represents a step to the left. 
\newline
Let {\it n} denote the total number of steps of the particle 
to the right until it is absorbed at the origin. (The number of 
steps to the left will then be $n+k$.)  
\newline
A path of length {\it n} --- or simply a path --- of the particle 
on this lattice is the set of horizontal and vertical bonds from 
the origin to the point of coordinates $(n, n+k)$ on the line 
{\it L}, as shown.
%\newline
%The variable {\it m} is not an independent
%variable, but depends on {\it n} and {\it k} by $m=n+k$. 
\newline
Let $C_{k}(n)$ be the {\it total} number of paths of length {\it n} 
from the origin to point $(n, n+k)$. A simple combinatorics analysis 
leads to the following expression for $C_{k}(n)$:
\be\label{fi}
C_{k}(n)=k~ \frac{(n+m-1)!}{{n!}~{m!}}=k~\frac{(2n+k-1)!}{{n!}(n+k)!}
=\frac{k}{(2n+k)} \left (
\begin{array}{cc}
2n+k \\
n
\end{array}
\right ).
\ee
The probability that 
%This number is related to the probabilities of absorption of the 
%particles at the origin. 
%A simple combinatorics analysis leads to the following expression for
the particle arrives at the origin given the initial
condition  that it starts at position $x=k$ is given by
\be\label{pf}
P(x=k)=\sum_{n=0}^{\infty}C_{k}(n) p^{n}(1-p)^{n+k}.
\ee
We show next how to carry out the summation in (\ref{pf}) using the
geometry of the lattice paths. We will need to solve explicitly for the
first three values of $k$, $k=1, 2, 3$, before arriving at the general 
solution (\ref{xi}) by induction.  
\vskip .3 cm
\noindent
\section{The case $k=1$} 
\begin{figure}
\begin{center}
\resizebox{!}{3.2 truecm}{\includegraphics{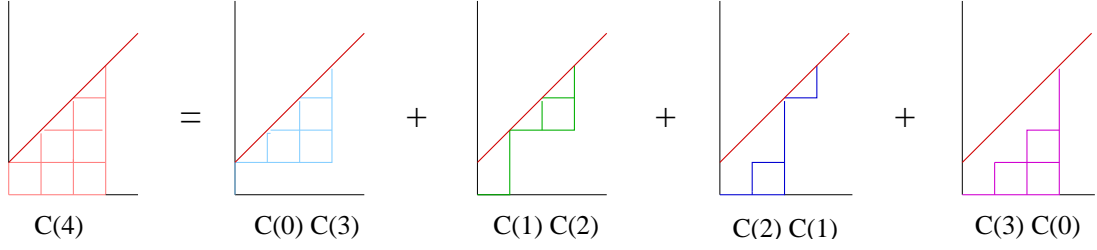}}
\vskip .3 cm
\parbox{13truecm}{\caption{\baselineskip=16 pt\small\label{fig2}
Decomposition of $C(n)$.}
}
\end{center}
\end{figure}

\vskip .3 cm
\noindent
The expression (\ref{fi}) for $C_{k=1}(n)=C(n)$ in the particular case 
the particle starts at $k=1$ is the $n^{th}$ Catalan number 
\be\label{11}
C(n)=\frac{(2n)!}{{n!}~{(n+1)!}}=\frac{1}{n+1} \left (
\begin{array}{cc}
2n \\
n
\end{array}
\right ).
\ee
The Catalan numbers have some geometric properties associated 
with paths on the lattice that we now explore. The first result
is the following
\vskip .1 cm
\noindent
{\theorem The number $C(n)$ given by (\ref{11}) satisfies the following 
equation 
\be\label{phin}
C(n)=\sum_{\alpha=1}^{n} C(\alpha-1)~C(n-\alpha).
\ee
}
\begin{itemize}
\item[] {\bf Proof.} Fig. \ref{fig2} illustrates the theorem for
$n=4$. For this proof, we have translated the paths in $C(4)$ one unit to
the left so as to make the paths touch the line {\it L}. These paths can
be decomposed as follows. Starting from $n=0$, and repeating the 
procedure for $n=1$, $n=2$ and $n=3$, raise vertical lines 
from the {\it n}-axis to the line {\it L}. The vertical line 
at $n=\alpha$ generates two sets of paths, the set 
$C(\alpha)$ to the left of $n=\alpha$, and the set 
$C(4-\alpha)$ to the right of $n=\alpha$. These two sets 
of paths contribute a factor $C(\alpha)C(4-\alpha)$ 
to the number $C(4)$. The sum of these factors for 
values from zero to $3$ is $C(4)$.
\hfill $\blacksquare$
\end{itemize} 
\subsection{Finding $P(x=1)$}
\vskip .3 cm
\noindent
Now if we set
\be\label{sum}
F(z)=\sum_{\alpha=1}^{\infty}~C(\alpha-1)~z^{\alpha},
\ee
then we obtain
\begin{eqnarray}\label{sum2}
F^{2}(z)&=&\sum_{\alpha,\beta=1}^{\infty}~C(\alpha-1)C(\beta-1)~z^{\alpha 
+\beta}  \nonumber\\
&=& 
\sum_{n=2}^{\infty}~\sum_{\alpha=1}^{n-1}~C(\alpha-1)C(n-\alpha 
-1)~z^{n} \nonumber\\
&=& \sum_{n=2}^{\infty}~C(n-1)~z^{n}=F(z)-z, \nonumber 
\end{eqnarray}
so that (\ref{sum}) obeys the equation 
\be\label{1}
F^{2}(z)=F(z)-z.
\ee
\noindent
Solving (\ref{1}) for {\it F} we get the solutions
\be
F(z)=\frac{1 \pm \sqrt{1-4z}}{2}.
\ee
Observe now that the summation (\ref{sum}) when $z=p(1-p)$ is (up to a
factor) the probability function (\ref{pf}) when $k=1$. Thus we  have
\be
P(x=1)=\frac{F(p-p^{2})}{p}=\frac{1 \pm
\sqrt{1-4p+4p^{2}}}{2p}=\frac{1 \pm (1-2p)}{2p}.
\ee
This last expression gives the classical results for the
probability that the particle will end up at the origin 
given that it starts at $k=1$. The solutions depend on
whether $p \leq  1/2$ or $p \geq 1/2$. For $p  \leq 1/2$, 
we get $P(x=1)=1$, and for $p  \geq 1/2$ we get 
\be\label{i}
P(x=1)=\frac{1-p}{p}.
\ee
For an interpretation of these results for other values of
{\it p}, see reference \cite{FM}.
\vskip .3 cm
\noindent
\section{The case $k=2$}
\begin{figure}
\begin{center}
\resizebox{!}{5 truecm}{\includegraphics{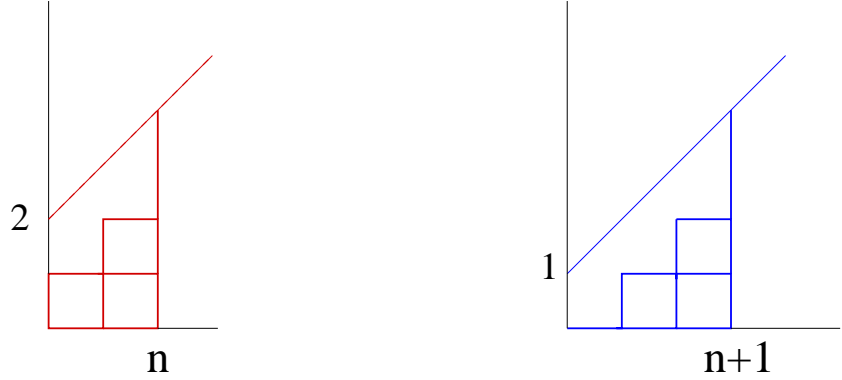}}
\vskip .3 cm
\parbox{13truecm}{\caption{\baselineskip=16 pt\small \label{fig3}
Relationship between $C_{2}(n)$ and the Catalan numbers.}
}
\end{center}
\end{figure}

\vskip .3 cm
\noindent
When $k=2$, then (\ref{pf}) becomes 
\be\label{pf1}
P(x=2)=\sum_{n=0}^{\infty}C_{2}(n) p^{n}(1-p)^{n+2}.
\ee
The coefficient $C_{2}(n)$ can be written in terms of $C(n)$ on
account of the following result
\vskip .3 cm
\noindent
{\theorem The following relationship holds for the coefficients $C_{2}(n)$
and $C(n)$ given by (\ref{fi}) 
\be\label{phin1}
C_{2}(n)=C(n+1).
\ee
}
\begin{itemize}
\item[] {\bf Proof.} In Fig. \ref{fig3} we have drawn the paths
that in $C_{2}(n)$ (in
red), and the paths in $C(n+1)$. Now, translate the paths in blue one 
unit to the left.  The blue paths go over the red ones. 
It is a perfect match and we are done.
\hfill $\blacksquare$
\end{itemize}
\vskip .5 cm
\noindent
If we now substitute (\ref{phin1}) into (\ref{pf1}) we obtain,
after some simple calculations,
\[
P(x=2)= \sum_{n=0}^{\infty} C(n+1)p^{n}(1-p)^{n+2} 
=\sum_{l=1}^{\infty} C(l)p^{l-1}(1-p)^{l+1}=
\frac{P(x=1)}{p}-\frac{1-p}{p} 
\]
and substituting for $P(x=1)$ from (\ref{i}) we get
\[
P(x=2)=\frac{(1-p)^{2}}{p^2} 
\]
if $p \geq 1/2$,  and $P(x=2)=1$ for $p \leq 1/2$.
\noindent
\vskip .3 cm
\noindent
\section{The case $k=3$}
\begin{figure}
\begin{center}
\resizebox{!}{5.4 truecm}{\includegraphics{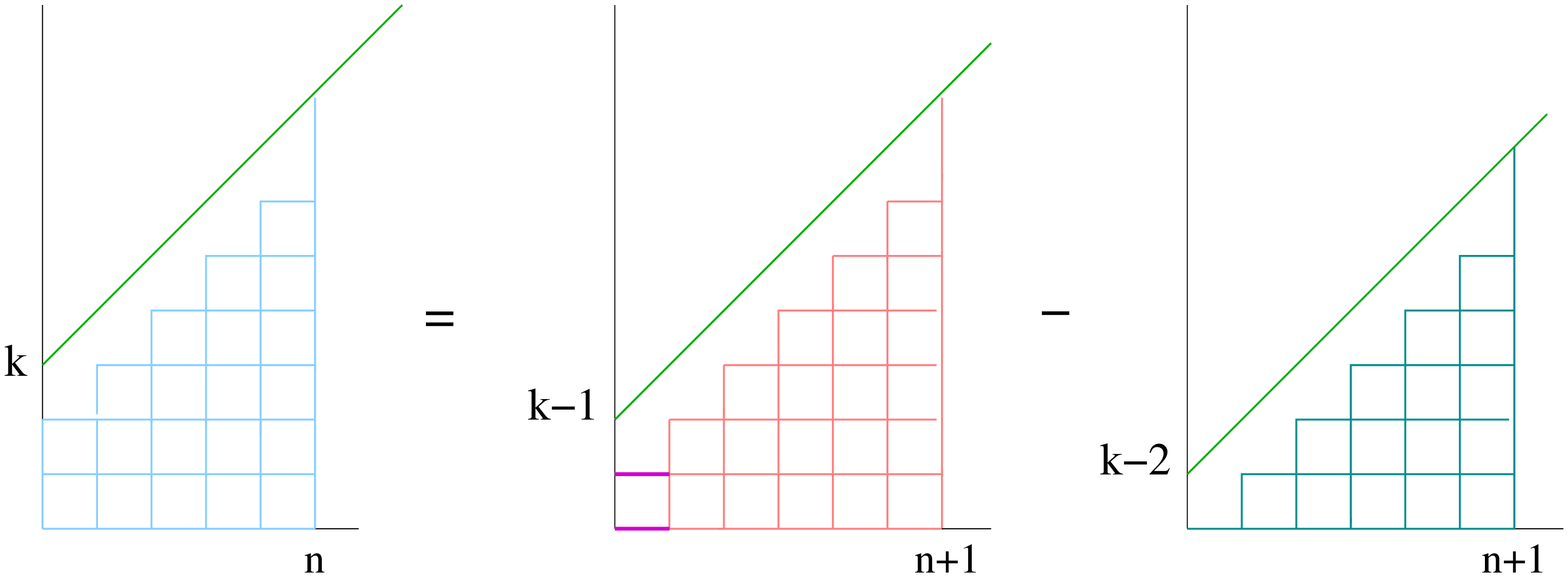}}
\vskip .3 cm
\parbox{13truecm}{\caption{\baselineskip=16 pt\small \label{fig4}
The recurrence relation (\ref{phin3}).}
}
\end{center}
\end{figure}

\vskip .3 cm
\noindent
We need this case to establish (\ref{xi}) for all values of {\it k}.  
\vskip .2 cm
\noindent
{\theorem
For $k \geq 3$ the following recurrence relation holds 
between the coefficient $C_{k}(n)$ and the coefficients $C_{k-1}(n+1)$ e
$C_{k-2}(n+1)$ 
\vskip .5 cm
\noindent
\be\label{phin3}
C_{k}(n)=C_{k-1}(n+1)-C_{k-2}(n+1).
\ee
}
\begin{itemize}
\item[] {\bf Proof.} The geometric representation of (\ref{phin3}) is
depicted in Fig. \ref{fig4}. To prove it, note that if we remove the paths
that stem from the two purple bonds on the right side of Fig. \ref{fig4} 
and translate the remaining set of paths one unit to the left, we get
exactly $C_{k}(n)$. But the number of paths we eliminated is exactly
$C_{k-2}(n+1)$, as the figure illustrates. This completes 
the proof.
\hfill $\blacksquare$
\end{itemize}
\vskip .5 cm
\noindent
The probability of absorption for $k=3$
\be\label{a}
P(x=3)=\sum_{n=0}^{\infty} C_{3}(n)p^{n}(1-p)^{n+3}.
\ee
Substituting for $C_{3}(n)$ its value from (\ref{phin3}),
\[
C_{3}(n)=C_{2}(n+1)-C_{1}(n+1), 
\]
we obtain, after some straightforward computation, 
\[
P(x=3)=\frac{(1-p)^{3}}{p^{3}}
\]
for $p \geq 1/2$ and $P(x=3)=1$ for $p \leq 1/2$. 
\vskip .3 cm
\noindent
\section{The general case} 
\vskip .3 cm
\noindent
By induction on $k$ we can now prove, as the three cases above suggest,
the following result.
\vskip .3 cm
\noindent
{\theorem The probability $P(x=k)$ of absorption of the particle at the
origin,
given that it starts at $x=k$, is given by 
\[
P(x=k)=\left ( P(x=1) \right )^{k}. 
\]
}
\begin{itemize}
\item[] {\bf Proof.} The probability of absorption of a particle that
starts at $x=k+1$ is equal to the probability that it moves one step to
the left to $x=k$ times the probability that it is absorbed at the 
origin from there, that is $(1-p)P(x=k)$, plus the probability 
that the particle moves one step to the right to $x=k+2$ times the
probability that it is absorbed from there, that is $ p P(x=k+2)$.
Hence, 
\[
P(x=k+1)=(1-p)P(x=k)+pP(x=k+2),
\]
from which we get
\be\label{p}
P(x=k+2)=\frac{1}{p}P(x=k+1)-\frac{(1-p)}{p}P(x=k),
\ee
The induction step now follows easily from (\ref{p}).
\hfill $\blacksquare$
\end{itemize}
\vskip .8 cm
\noindent
{\large \bf Acknowledgments\/}
I was supported by FAPESP -- Funda\c c\~ao de Amparo \`a Pesquisa do
Estado de S\~ao Paulo -- under grant 01/08485-6. I also thank 
Pierluigi Contucci for {\it other} path representations.

\end{document}